\documentclass[12pt]{article}

\usepackage{amsmath, epsfig, cite}
\usepackage{amssymb}
\usepackage{amsfonts}
\usepackage{latexsym}
\usepackage{float}
\usepackage{color}

\newtheorem{thm}{Theorem}[section]

\newtheorem{lem}[thm]{Lemma}

\newcommand*{\fl}[2]{\left\lfloor\frac{#1}{#2}\right\rfloor}

\numberwithin{equation}{section}

\newcommand{\qed}{{\hfill$\square$}\medskip}

\begin{document}

\begin{center}
{\Large\bf On a congruence involving $q$-Catalan numbers}
\end{center}

\vskip 2mm \centerline{Ji-Cai Liu}
\begin{center}
{\footnotesize Department of Mathematics, Wenzhou University, Wenzhou 325035, PR China\\
{\tt jcliu2016@gmail.com} }
\end{center}


\vskip 0.7cm \noindent{\bf Abstract.}
Based on a $q$-congruence of the author and Petrov, we set up a $q$-analogue of
Sun--Tauraso's congruence for sums of Catalan numbers, which extends a $q$-congruence due
to Tauraso.

\vskip 3mm \noindent {\it Keywords}: $q$-congruences; $q$-Catalan numbers; cyclotomic polynomials

\vskip 2mm
\noindent{\it MR Subject Classifications}: 11B65, 11A07, 05A10

\section{Introduction}
In combinatorics, the Catalan numbers are a sequence of natural numbers, which play an important role in various counting problems. The $n$th Catalan number is given by the following binomial coefficient:
\begin{align*}
C_n={2n\choose n}\frac{1}{n+1}={2n\choose n}-{2n\choose 2n+1}.
\end{align*}
The closely related numbers are the central binomial coefficients ${2n\choose n}$ for $n\ge 0$.

Both Catalan numbers and central binomial coefficients satisfy many interesting congruences (see, for instance, \cite{st-aam-2010,st-injt-2011}).
In 2011, Sun and Tauraso \cite{st-injt-2011} proved that for primes $p\ge 5$,
\begin{align}
&\sum_{k=0}^{p-1}{2k\choose k}\equiv \left(\frac{p}{3}\right)\pmod{p^2},\label{new-1}\\
&\sum_{k=0}^{p-1}C_k\equiv\frac{3}{2}\left(\frac{p}{3}\right)- \frac{1}{2}\pmod{p^2},\label{new-2}
\end{align}
where $\left(\frac{\cdot}{p}\right)$ denotes the Legendre symbol.

In the past few years, $q$-analogues of congruences ($q$-congruence) for indefinite sums of binomial coefficients as well as hypergeometric series attracted many experts' attention
(see, for example, \cite{gs-2019,guo-ijnt-2018,gz-aam-2010,gz-am-2019,lp-a-2019,Tauraso-aam-2012,tauraso-cm-2013}).
It is worth mentioning that Guo and Zudilin\cite{gz-am-2019} developed an interesting microscoping method to prove many $q$-congruences.

In order to discuss $q$-congruences, we first recall some $q$-series notation.
The $q$-binomial coefficients are defined as
\begin{align*}
{n\brack k}={n\brack k}_q
=\begin{cases}
\displaystyle\frac{(q;q)_n}{(q;q)_k(q;q)_{n-k}} &\text{if $0\leqslant k\leqslant n$},\\[10pt]
0 &\text{otherwise,}
\end{cases}
\end{align*}
where the $q$-shifted factorial is given by $(a;q)_n=(1-a)(1-aq)\cdots(1-aq^{n-1})$ for $n\ge 1$ and $(a;q)_0=1$. Moreover, the $q$-integers are defined by $[n]_q=(1-q^n)/(1-q)$, and the $n$th cyclotomic polynomial is given by
\begin{align*}
\Phi_n(q)=\prod_{\substack{1\le k \le n\\
(n,k)=1}}(q-e^{2k\pi i/n}).
\end{align*}

Recently, the author and Petrov \cite{lp-a-2019} established a $q$-analogue for \eqref{new-1} as follows:
\begin{align}
\sum_{k=0}^{n-1}q^k{2k\brack k}\equiv \left(\frac{n}{3}\right)q^{\frac{n^2-1}{3}}\pmod{\Phi_n(q)^2},\label{a-1}
\end{align}
which was originally conjectured by Guo \cite{guo-ijnt-2018} and generalises a $q$-congruence of
Tauraso \cite{Tauraso-aam-2012}. There are several natural $q$-analogues of Catalan numbers (see \cite{fh-jcta-1985}). Here and throughout the paper, we consider the following $q$-analogue of Catalan numbers:
\begin{align}
C_n(q)=\frac{1}{[n+1]_q}{2n\brack n}={2n\brack n}-q{2n\brack n+1}.\label{new-3}
\end{align}

In 2012, Tauraso\cite{Tauraso-aam-2012} obtained a weak $q$-version of \eqref{new-2} as follows:
\begin{align*}
\sum_{k=0}^{n-1}q^kC_k(q)\equiv \begin{cases}
\displaystyle q^{\lfloor n/3\rfloor}&\text{if $n\equiv 0,1 \pmod{3}$}\\
\displaystyle -1-q^{(2n-1)/3}&\text{if $n\equiv 2 \pmod{3}$}
\end{cases}\pmod{\Phi_n(q)},
\end{align*}
where $\lfloor x \rfloor$ denotes the integral part of real $x$.
In this note, we aim to set up a $q$-analogue of \eqref{new-2}
as well as another related $q$-congruence for sums of binomial coefficients.

\begin{thm}
For any positive integer $n$, the following holds modulo $\Phi_n(q)^2$:
\begin{align}
\sum_{k=0}^{n-1}q^{k}C_k(q)
\equiv
\begin{cases}
-q^{\frac{n^2-1}{3}}-q^{\frac{n(2n-1)}{3}}\quad &\text{if $n\equiv 2\pmod{3}$,}\\[7pt]
q^{\frac{n^2-1}{3}}-\frac{n-1}{3}(q^n-1)\quad &\text{if $n\equiv 1\pmod{3}$.}
\end{cases}\label{new-4}
\end{align}
\end{thm}

In order to prove \eqref{new-4}, we shall establish the following $q$-congruence.
\begin{thm}\label{t-1}
For any positive integer $n$, the following holds modulo $\Phi_n(q)^2$:
\begin{align}
\sum_{k=0}^{n-1}q^{k+1}{2k\brack k+1}
\equiv
\begin{cases}
q^{\frac{n(2n-1)}{3}}\quad &\text{if $n\equiv 2\pmod{3}$,}\\[7pt]
\frac{n-1}{3}(q^n-1)\quad &\text{if $n\equiv 1\pmod{3}$.}
\end{cases}\label{a-2}
\end{align}
\end{thm}

It is clear that \eqref{new-4} can be directly deduced from \eqref{a-1}, \eqref{new-3} and \eqref{a-2}. The remainder of the paper is organized as follows. We first set up a
preliminary result in the next section, and prove Theorem \ref{t-1} in Section 3.

\section{An auxiliary result}
\begin{lem}
For any positive integer $n$, the following holds modulo $\Phi_n(q)$:
\begin{align}
\sum_{k=1}^{n-1}\left(\frac{k-1}{3}\right)\frac{(-1)^{k}q^{\frac{1}{3}\left(2k^2-k
\left(\frac{k-1}{3}\right)\right)-\frac{k(k-1)}{2}}}{1-q^{k}}
\equiv
\begin{cases}
0\quad &\text{if $n\equiv 2\pmod{3}$},\\[7pt]
\frac{n-1}{6}\quad &\text{if $n\equiv 1\pmod{3}$}.
\end{cases}\label{b-1}
\end{align}
\end{lem}
{\it Proof.}
Note that
\begin{align*}
\sum_{k=1}^{n-1}(-1)^{k}\left(\frac{k-1}{3}\right)\frac{q^{\frac{1}{3}\left(2k^2-k
\left(\frac{k-1}{3}\right)\right)-\frac{k(k-1)}{2}}}{1-q^{k}}=
\sum_{k=0}^{\fl{n-3}{3}}\frac{(-1)^kq^{\frac{(k+1)(3k+2)}{2}}}{1-q^{3k+2}}
-\sum_{k=1}^{\fl{n-1}{3}}\frac{(-1)^kq^{\frac{k(3k+5)}{2}}}{1-q^{3k}}.
\end{align*}
We shall distinguish two cases to prove \eqref{b-1}.

{\noindent\bf Case 1}\quad $n\equiv 2\pmod{3}$.
This case is equivalent to
\begin{align}
\sum_{k=0}^{n-1}\frac{(-1)^kq^{\frac{(k+1)(3k+2)}{2}}}{1-q^{3k+2}}
-\sum_{k=1}^{n}\frac{(-1)^kq^{\frac{k(3k+5)}{2}}}{1-q^{3k}}
\equiv 0\pmod{\Phi_{3n+2}(q)}.\label{b-2}
\end{align}
Let $\omega$ be a primitive $(3n+2)$th root of unity. Letting $k\to n-k$ in the following sum gives
\begin{align*}
\sum_{k=0}^{n-1}\frac{(-1)^k\omega^{\frac{(k+1)(3k+2)}{2}}}{1-\omega^{3k+2}}
&=\sum_{k=1}^{n}\frac{(-1)^{n-k}\omega^{\frac{(n-k+1)(3n-3k+2)}{2}}}{1-\omega^{3n-3k+2}}\\
&=\sum_{k=1}^{n}\frac{(-1)^{n-k}\omega^{\frac{k(3k-1)}{2}+\frac{(3n+2)(n+1)}{2}-(3n+2)k}} {1-\omega^{3n-3k+2}}\\
&=\sum_{k=1}^{n}\frac{(-1)^{k}\omega^{\frac{k(3k+5)}{2}}}{1-\omega^{3k}},
\end{align*}
where we have used the fact that $\omega^{\frac{(3n+2)(n+1)}{2}}=(-1)^{n+1}$.
Thus,
\begin{align*}
\sum_{k=0}^{n-1}\frac{(-1)^k\omega^{\frac{(k+1)(3k+2)}{2}}}{1-\omega^{3k+2}}
-\sum_{k=1}^{n}\frac{(-1)^k\omega^{\frac{k(3k+5)}{2}}}{1-\omega^{3k}}
=0,
\end{align*}
which is equivalent to \eqref{b-2}.

{\noindent\bf Case 2}\quad $n\equiv 1\pmod{3}$. Let $\zeta$ be a primitive $(3n+1)$th root of unity. It suffices to show that
\begin{align}
\sum_{k=0}^{n-1}\frac{(-1)^k\zeta^{\frac{(k+1)(3k+2)}{2}}}{1-\zeta^{3k+2}}
-\sum_{k=1}^{n}\frac{(-1)^k\zeta^{\frac{k(3k+5)}{2}}}{1-\zeta^{3k}}=\frac{n}{2}.\label{b-3}
\end{align}
Note that
\begin{align*}
\sum_{k=0}^{n-1}\frac{(-1)^k\zeta^{\frac{(k+1)(3k+2)}{2}}}{1-\zeta^{3k+2}}
&=\sum_{k=n+1}^{2n}\frac{(-1)^{2n-k}\zeta^{\frac{(2n-k+1)(6n-3k+2)}{2}}}{1-\zeta^{6n-3k+2}}\\
&=\sum_{k=n+1}^{2n}\frac{(-1)^{k}\zeta^{\frac{k(3k-1)}{2}+(3n+1)(2n-2k+1)}}{1-\zeta^{-3k}}\\
&=-\sum_{k=n+1}^{2n}\frac{(-1)^{k}\zeta^{\frac{k(3k+5)}{2}}}{1-\zeta^{3k}},
\end{align*}
where we replace $k$ by $2n-k$ in the first step.
Thus,
\begin{align}
\sum_{k=0}^{n-1}\frac{(-1)^k\zeta^{\frac{(k+1)(3k+2)}{2}}}{1-\zeta^{3k+2}}
-\sum_{k=1}^{n}\frac{(-1)^k\zeta^{\frac{k(3k+5)}{2}}}{1-\zeta^{3k}}=
-\sum_{k=1}^{2n}\frac{(-1)^{k}\zeta^{\frac{k(3k+5)}{2}}}{1-\zeta^{3k}}.\label{b-4}
\end{align}
Furthermore, letting $k\to 2n+1-k$ on the right-hand side of \eqref{b-4} gives
\begin{align}
&\sum_{k=0}^{n-1}\frac{(-1)^k\zeta^{\frac{(k+1)(3k+2)}{2}}}{1-\zeta^{3k+2}}
-\sum_{k=1}^{n}\frac{(-1)^k\zeta^{\frac{k(3k+5)}{2}}}{1-\zeta^{3k}}\notag\\
&=-\sum_{k=1}^{2n}\frac{(-1)^{2n+1-k}\zeta^{\frac{(2n+1-k)(6n-3k+8)}{2}}}{1-\zeta^{3(2n+1-k)}}\notag\\
&=-\sum_{k=1}^{2n}\frac{(-1)^{1-k}\zeta^{\frac{(3k-1)(k-2)}{2}+(3n+1)(2n+3-2k)}}{1-\zeta^{1-3k}}\notag\\
&=-\sum_{k=1}^{2n}\frac{(-1)^{k}\zeta^{\frac{k(3k-1)}{2}}}{1-\zeta^{3k-1}}.\label{b-5}
\end{align}
An identity due to the author and Petrov \cite[(2.4)]{lp-a-2019} says
\begin{align}
\sum_{k=1}^{2n}\frac{(-1)^{k}\zeta^{\frac{k(3k-1)}{2}}}{1-\zeta^{3k-1}}=-\frac{n}{2}.\label{b-6}
\end{align}
Then the proof of \eqref{b-3} follows from \eqref{b-5} and \eqref{b-6}.
\qed

\section{Proof of Theorem \ref{t-1}}
Now we are in a position to prove Theorem \ref{t-1}. We recall the following identity:
\begin{align}
\sum_{k=0}^{n-1}q^k{2k\brack k+1}=\sum_{k=0}^{n-1}\left(\frac{n-k-1}{3}\right)q^{\frac{1}{3}\left(2(n-k)^2-(n-k)
\left(\frac{n-k-1}{3}\right)-3\right)}{2n\brack k},\label{c-1}
\end{align}
which was proved by Tauraso in a more general form (see \cite[Theorem 4.2]{Tauraso-aam-2012}).
Since $1-q^n\equiv 0\pmod{\Phi_n(q)}$, we have
\begin{align*}
1-q^{2n}=(1+q^n)(1-q^n)\equiv 2(1-q^n)\pmod{\Phi_n(q)^2}.
\end{align*}
It follows that for $1\le k\le n-1$,
\begin{align}
{2n\brack k}&=\frac{(1-q^{2n})(1-q^{2n-1})\cdots(1-q^{2n-k+1})}{(1-q)(1-q^2)\cdots(1-q^k)}\notag\\
&\equiv 2(1-q^{n})\frac{(1-q^{-1})\cdots(1-q^{-k+1})}{(1-q)(1-q^2)\cdots(1-q^k)}\pmod{\Phi_n(q)^2}\notag\\
&=2(q^n-1)\frac{(-1)^kq^{-\frac{k(k-1)}{2}}}{1-q^k}.\label{c-2}
\end{align}
Multiplying both sides of \eqref{c-1} by $q$ and substituting \eqref{c-2} into the right-hand side of \eqref{c-1}, we arrive at
\begin{align}
&\sum_{k=0}^{n-1}q^{k+1}{2k\brack k+1}\notag\\
&=\left(\frac{n-1}{3}\right)q^{\frac{1}{3}\left(2n^2-n
\left(\frac{n-1}{3}\right)\right)}+\sum_{k=1}^{n-1}\left(\frac{n-k-1}{3}\right)q^{\frac{1}{3}\left(2(n-k)^2-(n-k)
\left(\frac{n-k-1}{3}\right)\right)}{2n\brack k}\notag\\
&\equiv\left(\frac{n-1}{3}\right)q^{\frac{1}{3}\left(2n^2-n
\left(\frac{n-1}{3}\right)\right)}\notag\\
&+2(q^n-1)\sum_{k=1}^{n-1}\left(\frac{n-k-1}{3}\right)\frac{(-1)^kq^{\frac{1}{3}\left(2(n-k)^2-(n-k)
\left(\frac{n-k-1}{3}\right)\right)-\frac{k(k-1)}{2}}}{1-q^k}\pmod{\Phi_n(q)^2}.\label{c-3}
\end{align}

Furthermore,
\begin{align*}
&\sum_{k=1}^{n-1}\left(\frac{n-k-1}{3}\right)\frac{(-1)^kq^{\frac{1}{3}\left(2(n-k)^2-(n-k)
\left(\frac{n-k-1}{3}\right)\right)-\frac{k(k-1)}{2}}}{1-q^k}\notag\\
&=\sum_{k=1}^{n-1}\left(\frac{k-1}{3}\right)\frac{(-1)^{n-k}q^{\frac{1}{3}\left(2k^2-k
\left(\frac{k-1}{3}\right)\right)-\frac{(n-k)(n-k-1)}{2}}}{1-q^{n-k}}\\
&=\sum_{k=1}^{n-1}\left(\frac{k-1}{3}\right)\frac{(-1)^{n-k}q^{\frac{1}{3}\left(2k^2-k
\left(\frac{k-1}{3}\right)\right)-\frac{n(n-1)}{2}-\frac{k(k+1)}{2}+nk}}{1-q^{n-k}}\\
&\equiv \sum_{k=1}^{n-1}\left(\frac{k-1}{3}\right)\frac{(-1)^{k}q^{\frac{1}{3}\left(2k^2-k
\left(\frac{k-1}{3}\right)\right)-\frac{k(k-1)}{2}}}{1-q^{k}}\pmod{\Phi_n(q)},
\end{align*}
where we set $k\to n-k$ in the first step.
Thus,
\begin{align}
\sum_{k=0}^{n-1}q^{k+1}{2k\brack k+1}&\equiv
\left(\frac{n-1}{3}\right)q^{\frac{1}{3}\left(2n^2-n
\left(\frac{n-1}{3}\right)\right)}\notag\\
&+2(q^n-1)\sum_{k=1}^{n-1}\left(\frac{k-1}{3}\right)\frac{(-1)^{k}q^{\frac{1}{3}\left(2k^2-k
\left(\frac{k-1}{3}\right)\right)-\frac{k(k-1)}{2}}}{1-q^{k}}\pmod{\Phi_n(q)^2}.\label{new-5}
\end{align}
We complete the proof of \eqref{a-2} by combining \eqref{b-1} and \eqref{new-5}.

\vskip 5mm \noindent{\bf Acknowledgments.}
This work was supported by the National Natural Science Foundation of China (grant 11801417).

\end{document}